# OPTIMAL DESIGNS FOR THREE-DIMENSIONAL SHAPE ANALYSIS WITH SPHERICAL HARMONIC DESCRIPTORS[1]

By Holger Dette, Viatcheslav B. Melas[2]
and Andrey Pepelyshev

*Ruhr-Universität Bochum, St. Petersburg State University
and St. Petersburg State University*

We determine optimal designs for some regression models which are frequently used for describing three-dimensional shapes. These models are based on a Fourier expansion of a function defined on the unit sphere in terms of spherical harmonic basis functions. In particular, it is demonstrated that the uniform distribution on the sphere is optimal with respect to all $\Phi_p$ criteria proposed by Kiefer in 1974 and also optimal with respect to a criterion which maximizes a $p$ mean of the $r$ smallest eigenvalues of the variance–covariance matrix. This criterion is related to principal component analysis, which is the common tool for analyzing this type of image data. Moreover, discrete designs on the sphere are derived, which yield the same information matrix in the spherical harmonic regression model as the uniform distribution and are therefore directly implementable in practice. It is demonstrated that the new designs are substantially more efficient than the commonly used designs in three-dimensional shape analysis.

**1. Introduction.** Over the last decade, tools for acquiring and visualizing three-dimensional (3D) models have become integral components of data processing in many fields, including medicine, chemistry, architecture, agriculture and biology. Volumetric shape analysis permits an evaluation of the actual structures that are implicitly represented in 3D image data. For the analysis, description and comparison of shapes of various structures, shape descriptors, which are able to handle very different shapes and

Received March 2004; revised January 2005.

[1]Supported by the Deutsche Forschungsgemeinschaft SFB 475, "Komplexitätsreduktion in multivariaten Datenstrukturen."

[2]Supported in part by the Russian Foundation of Basic Reseach Project 04-01-00519.

---

to represent their global and local features, are of increasing interest (see, e.g., Brechbühler, Gerig and Kübler [4], Novotni and Klein [14], Székely, Kelemen, Brechbühler and Gerig [23], Ding, Nesumi, Takano and Ukai [5], Funkhouser, Min, Kazhdan, Chen, Halderman, Dobkin and Jacobs [7] and Kazhdan, Funkhouser and Rusinkiewicz [11], among many others). Spherical harmonic shape descriptors usually describe the surface in terms of a relatively small number of coefficients of a spherical harmonic expansion of the radius as a function on the unit sphere (see, e.g., [5] or [4]), that is,

$$(1.1) \qquad r(\theta, \phi) = \sum_{\ell=0}^{\infty} \sum_{m=-\ell}^{\ell} c_\ell^m Y_\ell^m(\theta, \phi),$$

where $\theta \in [0, \pi], \phi \in (-\pi, \pi]$, the quantities

$$(1.2) \qquad c_\ell^m = \frac{1}{4\pi} \int_0^\pi \int_{-\pi}^\pi r(\theta, \phi) Y_\ell^m(\theta, \phi) \, d\phi \sin\theta \, d\theta$$

are the usual "Fourier" coefficients and

$$\{Y_\ell^m(\theta, \phi) | m \in \{-\ell, -\ell+1, \dots, \ell\}; \ell \in \mathbb{N}_0\}$$

is a complete orthonormal basis on the unit sphere. Let $r_i = r(\theta_i, \phi_i)$ denote the observed radius of the 3D shape at polar angle $\theta_i$ and azimuthal angle $\phi_i$ [in other words, the corresponding point of the shape has spherical coordinates $(r_i \sin\theta_i \cos\phi_i, r_i \sin\theta_i \sin\phi_i, r_i \cos\theta_i)^T$] and assume that data

$$\{(r_i, \theta_i, \phi_i) | i = 1, \dots, n\}$$

are available for one object. Usually a truncated expansion of order $d$ is applied as an approximation of (1.1), where the coefficients $c_\ell^m$ are determined by the least squares critierion

$$(1.3) \qquad \min_{c_\ell^m} \left\{ \sum_{i=1}^n \left( r_i - \sum_{\ell=0}^d \sum_{m=-\ell}^\ell c_\ell^m Y_\ell^m(\theta_i, \phi_i) \right)^2 \right\}$$

and the estimated coefficients in this expansion (appropriately normalized) are then used for describing and analyzing the 3D shapes. For this purpose typical tools from multivariate statistics (cluster, discriminant or principal component analysis) are applied (see Ding, Nesumi, Takano and Ukai [5], Kazhdan, Funkhouser and Rusinkiewicz [11] and Kelemen, Székely and Gerig [12], among many others).

A common experimental design for this situation is a uniform distribution on the rectangle $[0, \pi] \times (-\pi, \pi]$ realized by a grid or a uniform design that takes observations on several circles with equal distance on the $z$ axis (see, e.g., [5]). In the literature on shape analysis these designs are mainly motivated by their easy implementation. If the grid is fine enough or a sufficiently large number of circles on the unit sphere are used, the matrix $B^T B$ of the least squares estimate $\hat{c} = (B^T B)^{-1} B^T r$ approximates a diagonal matrix, which simplifies the numerical calculation in the statistical analysis. Here $r = (r_1, \ldots, r_n)$ is the vector of measured radii and

$$B = (Y_\ell^m(\theta_i, \phi_i))_i^{\ell, m}$$

is the design matrix corresponding to the least squares problem (1.3) (see [4]).

In the present paper we consider the problem of finding optimal designs for 3D shape analysis based on spherical harmonic descriptors. In Section 2 we present some more details on spherical harmonic descriptors and basic results on the theory of optimal experimental design. In the same section we also demonstrate that the uniform distribution on the unit sphere is optimal with respect to any of Kiefer's [13] $\Phi_p$ criteria if the interest of the experimenter is the estimation of the complete vector

$$c = (c_0^0, c_1^{-1}, c_1^0, c_1^1, \ldots, c_d^{-d})^T \in \mathbb{R}^{(d+1)^2}$$

or a certain subset of the parameters. It is also shown that for any $t \leq (d+1)^2$ this design maximizes the $p$th mean of the $r$ largest eigenvalues of the variance–covariance matrix $(B^T B)^{-1}$. Therefore the uniform distribution on the sphere is particularly efficient for principal component analysis, which is the main tool for summarizing the information contained in the spherical harmonic coefficients obtained by (1.3) (see, e.g., [5] or [11]). Because this design is continuous with density $\frac{1}{4\pi} \sin \theta \, d\theta \, d\phi$, it is not directly implementable in practice. Therefore, for a finite sample size we determine in Section 3 discrete designs which give the same information matrix as the uniform distribution on the sphere. For this reason these designs are also optimal with respect to Kiefer's [13] $\Phi_p$ criteria, and optimal with respect to a criterion that maximizes a $p$ mean of the $r$ smallest eigenvalues of the information matrix and is related to principal component analysis. In Section 4 we present several examples which illustrate the advantages of our approach, and determine optimal uniform designs which take in each direction only one observation and are for this reason particularly attractive to practitioners. We also reanalyze a design used by Ding, Nesumi, Takano and Ukai [5] for principal component analysis and demonstrate that the designs derived in the present paper are substantially more efficient. Finally, some conclusions and directions for future research are mentioned in Section 5, while more technical details are given in the Appendix.

## 2. Spherical harmonic descriptors and optimal design.

An orthogonal system

$$\{Y_\ell^m(\theta, \phi) | \ell \in \mathbb{N}_0; m \in \{-\ell, -\ell+1, \ldots, \ell\}\}$$

of functions on the unit sphere satisfies

$$(2.1) \qquad \frac{1}{4\pi} \int_0^\pi \int_{-\pi}^\pi Y_\ell^m(\theta,\phi) Y_{\ell'}^{m'}(\theta,\phi) \, d\phi \sin\theta \, d\theta = 0$$

whenever $(m,\ell) \neq (m',\ell')$. The common system used in shape analysis (see [4] or [5]) is obtained by the normalization

$$\frac{1}{4\pi} \int_0^\pi \int_{-\pi}^\pi (Y_\ell^m(\theta,\phi))^2 \, d\phi \sin\theta \, d\theta = 1$$

and given by the spherical harmonic functions

$$Y_n^0(\theta,\phi) = \sqrt{2n+1} P_n^0(\cos\theta), \qquad n \in \mathbb{N}_0,$$

$$Y_n^m(\theta,\phi) = \sqrt{2(2n+1)\frac{(n-m)!}{(n+m)!}} P_n^m(\cos\theta) \cos(m\phi),$$

$$(2.2) \qquad\qquad\qquad\qquad\qquad\qquad m = 0,\ldots,n; n \in \mathbb{N},$$

$$Y_n^{-m}(\theta,\phi) = \sqrt{2(2n+1)\frac{(n+m)!}{(n-m)!}} P_n^{-m}(\cos\theta) \sin(m\phi),$$

$$\qquad\qquad\qquad\qquad\qquad\qquad m = -n,\ldots,-1; n \in \mathbb{N},$$

where $P_n^m(x)$ is the $m$th associated Legendre function of degree $n$ satisfying the differential equation

$$(2.3) \qquad (1-x^2)P''(x) - 2xP'(x) + \left\{ n(n+1) - \frac{m^2}{1-x^2} \right\} P(x) = 0.$$

It is well known (see [1], Chapter 9) that these functions can be represented as

$$(2.4) \qquad P_n^m(x) = (-1)^m \frac{(2m)!}{2^n m!} (1-x^2)^{m/2} C_{n-m}^{(m+1/2)}(x),$$

where $C_k^{(\alpha)}(x)$ is the $k$th ultraspherical polynomial orthogonal with respect to the measure $(1-x^2)^{\alpha-1/2} \, dx$ (see [22]). Note also that $P_n^0(x)$ is the $n$th Legendre polynomial $P_n(x)$ orthogonal with respect to the Lebesgue measure on the interval $[-1,1]$. Because orthogonal polynomials can be calculated recursively, this representation allows a fast computation of the functions $Y_\ell^m$, and the first four spherical harmonic functions corresponding to the case $d=1$ are given by

$$(2.5) \qquad \begin{aligned} &Y_0^0(\theta,\phi) = 1, & &Y_1^0(\theta,\phi) = \sqrt{3}\cos\theta, \\ &Y_1^{-1}(\theta,\phi) = \sqrt{3}\sin\theta\sin\phi, & &Y_1^1(\theta,\phi) = \sqrt{3}\sin\theta\cos\phi. \end{aligned}$$

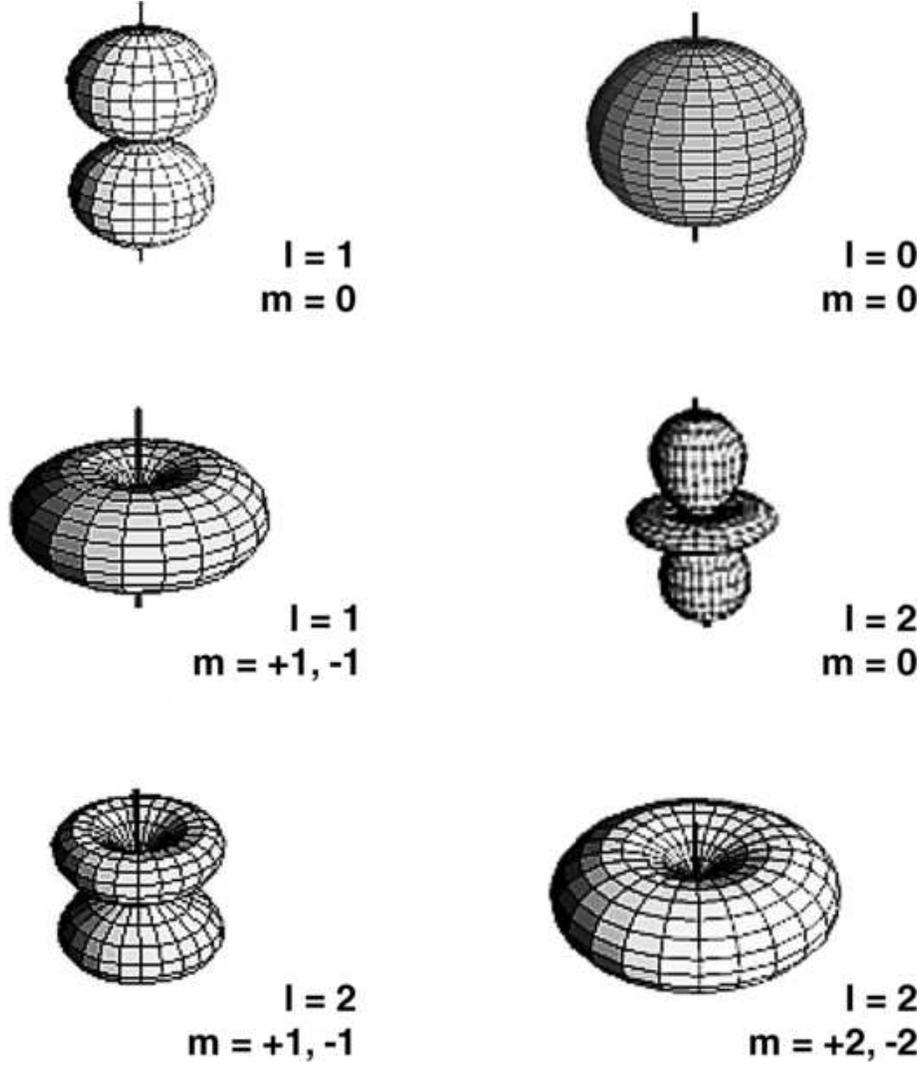

Fig. 1. *Spherical harmonic descriptors* $(Y_\ell^m(\theta, \phi) \sin \theta \cos \phi, Y_\ell^m(\theta, \phi) \sin \theta \sin \phi, Y_\ell^m(\theta, \phi) \times \cos \theta)^T$ *for* $m = -\ell, \ldots, \ell,\ \ell = 0, 1, 2.$

Figures 1 and 2 show the spherical harmonic descriptors

$$(Y_\ell^m(\theta, \phi) \sin \theta \cos \phi, Y_\ell^m(\theta, \phi) \sin \theta \sin \phi, Y_\ell^m(\theta, \phi) \cos \theta)^T, \qquad m = -\ell, \ldots, \ell,$$

for $\ell = 0, 1, 2, 3$, when $(\theta, \phi)$ varies in the rectangle $[0, \pi] \times (-\pi, \pi]$.

Consider the regression model corresponding to the least squares problem (1.3),

$$(2.6) \qquad E[Y|\theta, \phi] = c^T f_d(\theta, \phi), \qquad \text{Var}[Y|\theta, \phi] = \sigma^2 > 0,$$

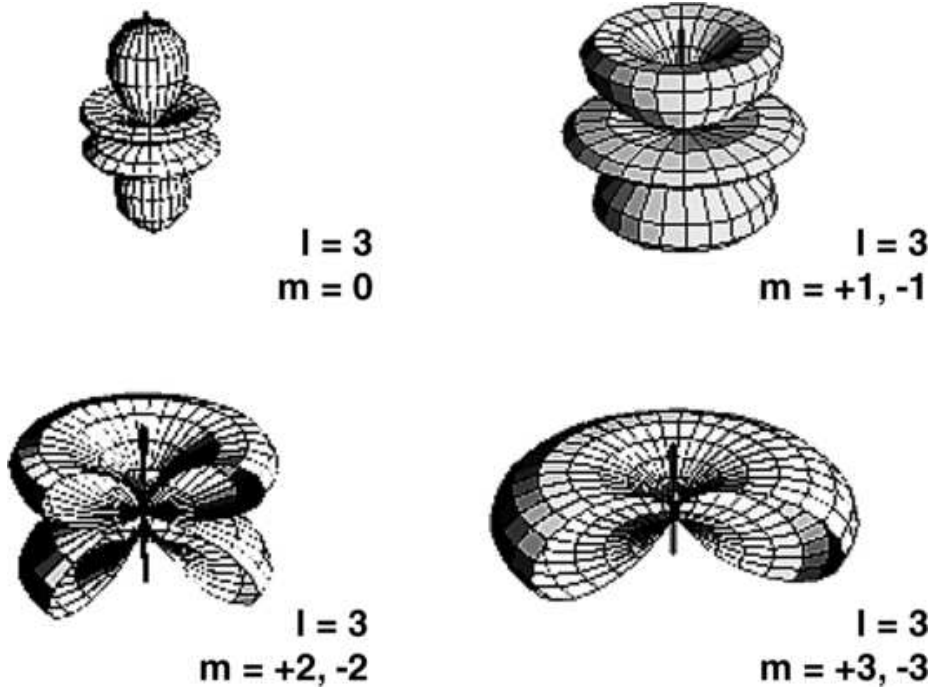

FIG. 2. *Spherical harmonic descriptors* $(Y_\ell^m(\theta,\phi)\sin\theta\cos\phi, Y_\ell^m(\theta,\phi)\sin\theta\sin\phi, Y_\ell^m(\theta,\phi)\times\cos\theta)^T$ *for* $m=-\ell,\ldots,\ell$, $\ell=3$.

where

$$
\begin{aligned}
f_d(\theta,\phi) = (Y_0^0(\theta,\phi), Y_1^{-1}(\theta,\phi), Y_1^0(\theta,\phi),\\
Y_1^1(\theta,\phi),\ldots,Y_d^{-d}(\theta,\phi),\ldots,Y_d^d(\theta,\phi))^T \in \mathbb{R}^{(d+1)^2}
\end{aligned}
\tag{2.7}
$$

is the vector of spherical harmonic functions of order $d$ and

$$
c = (c_0^0, c_1^{-1}, c_1^0, c_1^1,\ldots,c_d^{-d},\ldots,c_d^d)^T \in \mathbb{R}^{(d+1)^2}
$$

is the corresponding vector of parameters. Note that $(d+1)^2$ spherical harmonic functions appear in the model (2.6). An approximate design is a probability measure on the set $[0,\pi]\times(-\pi,\pi]$. For a probability measure with finite support, the support points, say $z_i = (\theta_i,\phi_i)$, determine the points on the sphere where the radius of the 3D shape is observed, and the corresponding weights, say $w_i$, give the relative proportion of total observations taken in a particular direction. For a given design $\xi$ with finite support, the covariance matrix of the least squares estimate for the vector $c$ is approximately proportional to the information matrix

$$
M(\xi) = \int_{-\pi}^{\pi}\int_0^{\pi} f(\theta,\phi) f^T(\theta,\phi)\,d\xi(\theta,\phi)
\tag{2.8}
$$

(this is essentially the matrix $B^T B$ mentioned in the Introduction), and an optimal approximate design maximizes an appropriate function of this matrix. There are numerous criteria proposed in the literature that can be used to discriminate between competing designs (see [19] or [16]), and we will restrict ourselves to the famous family of $\Phi_p$ criteria introduced by Kiefer [13], and a new optimality criterion which has not been considered so far in the literature and is related to principal component analysis.

To be precise, let $K \in \mathbb{R}^{(d+1)^2 \times s}$ denote a given matrix of rank $s \leq (d+1)^2$ and assume that the main interest of the experimenter is in the estimation of $s$ linear combinations $K^T c$. If $n$ observations are taken according to an approximate design (possibly by applying an appropriate rounding procedure to a design with finite support or by using a discretization of a continuous design; see [17]), then the covariance matrix of the least squares estimate for $K^T c$ is approximately given by

$$\frac{\sigma^2}{n}(K^T M^-(\xi)K),$$

where $A^-$ denotes a generalized inverse of the matrix $A$ and we assume that the linear combinations $K^T c$ are estimable by the design $\xi$, that is,

$$\text{range}(K) \subset \text{range}(M(\xi));$$

see [16]. Let $-\infty \leq p < 1$. Following [13], we call the design $\xi^*$ $\Phi_p$-optimal for estimating the linear combinations $K^T c$ if $\xi^*$ maximizes the expression

$$(2.9) \qquad \Phi_p(\xi) = (\text{tr}(K^T M^-(\xi)K)^{-p})^{1/p}$$

among all designs for which $K^T c$ is estimable. If $K = I_{(d+1)^2}$ is the identity matrix of order $(d+1)^2 \times (d+1)^2$, then $\xi^*$ is briefly called $\Phi_p$-optimal. Note that the cases $p = 0$ and $p = -\infty$ correspond to the frequently used $D$- and $E$-optimality criteria, that is,

$$(2.10) \quad \Phi_0(\xi) = \det(K^T M^-(\xi)K)^{-1}, \qquad \Phi_{-\infty}(\xi) = \lambda_{\min}((K^T M^-(\xi)K)^{-1}),$$

while the $A$ criterion is obtained for the choice $p = -1$, that is,

$$(2.11) \qquad \Phi_{-1}(\xi) = \text{tr}(K^T M^-(\xi)K)^{-1}.$$

In the following discussion, we are interested in a design that is particularly efficient for the estimation of the coefficients corresponding to the $(2k + 1)$ spherical harmonic functions

$$Y_k^{-k}, Y_k^{-k+1}, \ldots, Y_k^k$$

of the vector of regression functions defined in (2.7), where $k \in \{0, \ldots, d\}$ denotes a given "level of resolution." For this, define $0_{k,s}$ as the $(2k + 1) \times$ $(2s + 1)$ matrix with all entries equal to $0$, let $q \in \mathbb{N}_0$, consider the indices $0 \le k_0 < k_1 < k_2 < \cdots < k_q \le d$ and define the matrix

$$(2.12) \qquad K^T = K^T_{k_0, \ldots, k_q} = (K_{i,j})^{j=0,\ldots,d}_{i=0,\ldots,q}$$

by

$$(2.13) \qquad K_{i,j} = \begin{cases} 0_{k_i, j}, & \text{if } j \neq k_i, \\ I_{2k_i + 1}, & \text{if } j = k, \end{cases}$$

$(i = 0, \ldots, q; j = 0, \ldots, d)$. Note that $K \in \mathbb{R}^{(d+1)^2 \times s}$, where $s = \sum_{i=0}^q (2k_i + 1)$, and that $K^T c$ gives the vector with coefficients

$$(2.14) \qquad \{c^m_{k_\ell} | m \in \{-k_\ell, -k_\ell + 1, \ldots, k_\ell\}; \ell = 0, \ldots, q\}.$$

Two cases are of particular interest here and are therefore mentioned separately. If $q = d$, then

$$(2.15) \qquad K^T_{0, \ldots, d} = I_{(d+1)^2}$$

and precise estimation of the full vector of parameters is the main goal for the construction of the optimal design. On the other hand, if $q = 0$, only the coefficients corresponding to the $(2k_0 + 1)$ spherical harmonic functions $Y^{-k_0}_{k_0}, \ldots, Y^{k_0}_{k_0}$ are of interest and the corresponding matrix is given by

$$(2.16) \quad K^T_{k_0} = [0_{k_0, 0} \vdots \cdots \vdots 0_{k_0, k_0 - 1} \vdots I_{2k_0 + 1} \vdots \cdots \vdots 0_{k_0, d}] \in \mathbb{R}^{(2k_0 + 1) \times (d+1)^2}.$$

Note that the general matrix defined by (2.12) consists of $(q + 1)$ (block) rows of the form (2.16). The following result shows that for matrices of this form, the $\Phi_p$ criteria defined in (2.9) are maximized by the uniform distribution on the sphere independently of $p \in [-\infty, 1)$.

THEOREM 2.1. *Let $p \in [-\infty, 1)$, let $0 \le k_0 < \cdots < k_q \le d$ be given indices and denote $K = K_{k_0, \ldots, k_q}$ as the matrix defined by (2.12) and (2.13). A $\Phi_p$-optimal design $\xi^*$ for estimating the linear combinations $K^T c$ in the spherical harmonic regression model (2.6) is given by the uniform distribution on the sphere, that is,*

$$(2.17) \qquad \xi^*(d\theta, d\phi) = \frac{1}{4\pi} \sin\theta \, d\theta \, d\phi.$$

*Moreover, the corresponding information matrix in the spherical harmonic regression model is given by $M(\xi^*) = I_{(2d+1)^2}$.*

PROOF. Let $\xi^*$ denote the design corresponding to the density (2.17). Then, due to the orthonormality of the spherical harmonic functions $Y^m_\ell$, it follows that

$$M(\xi^*) = \left( \int_0^\pi \int_{-\pi}^\pi Y^m_\ell(\theta, \phi) Y^{m'}_{\ell'}(\theta, \phi) \, d\phi \sin\theta \, \frac{d\theta}{4\pi} \right)_{\ell, \ell', m', m} = I_{(2d+1)^2}.$$

Assume for the moment that $p > -\infty$. According to Theorem 7.20 in [16], we obtain that the measure $\xi^*$ is $\Phi_p$-optimal if and only if the inequality

$$(2.18) \qquad f_d^T(\theta, \phi) K (K^T K)^{-p-1} K^T f(\theta, \phi) \leq \operatorname{tr}(K^T K)^{-p}$$

holds for all $\theta \in [0, \pi]$ and $\phi \in (-\pi, \pi]$. Now the special structure of the matrix $K$ defined in (2.12) and (2.13) yields

$$K^T K = I_s,$$

where $s = \sum_{\ell=0}^{q} (2k_\ell + 1)$. Observing the definition of the vector $f_d(\theta, \phi)$ in (2.7) and again the definition of the matrix $K$, (2.18) reduces to

$$
\begin{aligned}
(2.19) \qquad s &\geq \sum_{\ell=0}^{q} \sum_{m=-k_\ell}^{k_\ell} (Y_{k_\ell}^m(\theta, \phi))^2 \\
&= \sum_{\ell=0}^{q} (2k_\ell + 1) \left\{ (P_\ell^0(\cos\theta))^2 + 2 \sum_{m=1}^{k_\ell} \frac{(k_\ell - m)!}{(k_\ell + m)!} (P_{k_\ell}^m(\cos\theta))^2 \right\}
\end{aligned}
$$

for all $\theta \in [0, \pi]$ and $\phi \in (-\pi, \pi]$, where we have used the representation (2.2) and the trigonometric identity $\cos^2(m\phi) + \sin^2(m\phi) = 1$ for the last equality. From the identity

$$P_k^0(\cos\alpha\cos\beta + \sin\alpha\sin\beta\cos\phi)$$

$$= P_k^0(\cos\alpha) P_k^0(\cos\beta) + 2 \sum_{m=1}^{k} \frac{(k-m)!}{(k+m)!} P_k^m(\cos\alpha) P_k^m(\cos\beta) \cos(m\phi)$$

for the Legendre functions (see [1], page 457) and the fact that $P_k^0$ is the $k$th Legendre polynomial, it follows (with the choice $\alpha = \beta, \phi = 0$) that

$$(P_k^0(\cos\alpha))^2 + 2 \sum_{m=1}^{k} \frac{(k-m)!}{(k+m)!} (P_k^m(\cos\alpha))^2 = P_k^0(1) = 1.$$

Now the right-hand side of (2.19) simplifies to

$$\sum_{\ell=0}^{q} (2k_\ell + 1) = s$$

and, consequently, (2.19) or equivalently (2.18) holds for all $\theta \in [0, \pi]$ and $\phi \in (-\pi, \pi]$. This proves that the design $\xi^*$ corresponding to the uniform distribution on the sphere is $\Phi_p$-optimal for any $p \in (-\infty, 1)$ and any matrix $K$ of the form (2.12) and (2.13). The remaining case $p = -\infty$ finally follows from Lemma 8.15 in [16]. $\square$

It should be noted at this point that the $E$ criterion is of particular importance for the design of experiments in 3D shape analysis, because many authors propose to use principal component analysis for the comparison of different shapes (see, e.g., Ding, Nesumi, Takano and Ukai [5] or Kazhdan, Funkhouser and Rusinkiewicz [11], among others). More precisely, it is proposed to summarize the information contained in the spherical harmonic coefficients by a principal component analysis based on the variance–covariance matrix of the least squares estimator obtained from (1.3), which explains the 3D shape variation of different objects by the first, say $r$, principal components corresponding to the $r$ largest eigenvalues of the matrix $M^{-1}(\xi)$. Consequently, an efficient design for principal component analysis should minimize a function of the largest $r$ eigenvalues of the matrix $M^{-1}(\xi)$. To be precise, let $r \leq (d+1)^2$ and $-p \leq \infty < 1$. Then we call a design $\xi^*$ $\Psi_{p,r}$-optimal if $\xi^*$ maximizes

$$(2.20) \qquad \Psi_{p,r}(\xi) = \left( \sum_{j=1}^{r} \{ \lambda_{(j)}(M(\xi)) \}^p \right)^{1/p},$$

where $\lambda_{(j)}(M(\xi))$ denotes the $j$th smallest eigenvalue of the matrix $M(\xi^*)$. Again the cases $p = -\infty$ and $p = 0$ are obtained from the corresponding limits, that is, $\Psi_{-\infty,r}(\xi) = \lambda_{(1)}(M(\xi))$ and $\Psi_{0,r}(\xi) = \prod_{j=1}^{r} \lambda_{(j)}(M(\xi))$, respectively. To our knowledge, optimality criteria of this type have not been considered in the literature so far. Note also that the $E$- and $A$-optimality criteria (with $K = I_{(d+1)^2}$) are obtained for the choices $r = 1$ and $r = (d+1)^2$ with $p = -1$, respectively. The most important case for the choice of the parameter $p \in [-\infty, 1)$ is certainly obtained for $p = -1$. In the following corollary, we show that the uniform distribution on the sphere is also $\Psi_{r,p}$-optimal.

COROLLARY 2.2. *Let $p \in [-\infty, 1)$ and $1 \leq r \leq (d+1)^2$. Then the uniform distribution on the sphere defined by* (2.17) *is $\Psi_{p,r}$-optimal in the spherical harmonic regression model* (2.6).

PROOF. The case $p = -\infty$ is obtained from Theorem 2.1 with $K = I_{(d+1)^2}$. For $-\infty < p < 1$, note that $\Psi_{p,r}(\xi^*) = r^{1/p}$ if $p \neq 0$ and $\Psi_{0,r}(\xi^*) = 1$. If $p \neq 0$ and $\xi^*$ were not $\Psi_{p,r}$-optimal, then there would exist a design, say $\xi$, with

$$(\Psi_{p,r}(\xi))^p = \sum_{j=1}^{r} \lambda_{(j)}^p(M(\xi)) < r = (\Psi_{p,r}(\xi^*))^p.$$

Consequently, we obtain $\lambda_{(1)}(M(\xi)) < 1 = \lambda_{\min}(M(\xi^*))$, which implies that the uniform distribution on the sphere $\xi^*$ is not $E$-optimal and contradicts Theorem 2.1. The case $p = 0$ is obtained by a similar argument and the assertion of the corollary has been established. $\square$

**3. Discrete optimal designs.** Note that the uniform distribution defined by (2.17) is not directly implementable as a design in real experiments. Therefore, for practical applications it is important to obtain discrete designs $\xi$ which are equivalent to the uniform distribution $\xi^*(d\theta, d\phi) = \frac{1}{4\pi} \sin\theta \, d\theta \, d\phi$ in the sense that

$$(3.1) \qquad M(\xi) = M(\xi^*) = I_{(d+1)^2},$$

where $M(\xi)$ is the information matrix (2.8) of the design $\xi$ in the spherical harmonic regression model (2.6). Note that due to Carathéodory's theorem (see [19]), there always exist discrete designs satisfying (3.1), and in the following we will construct a broad class of discrete designs that can easily be implemented in practice. For this purpose we need an auxiliary result about quadrature formulas, which will be proved in the Appendix.

LEMMA 3.1. *Assume that $r \in \mathbb{N}$ and let $-1 \leq x_1 < x_2 < \cdots < x_r \leq 1$ denote $r$ points with corresponding weights $w_1, \ldots, w_r > 0$ ($\sum_{i=1}^r w_i = 1$). The points $x_i$ and weights $w_i$ generate a quadrature formula of degree $z \geq r$, that is,*

$$(3.2) \qquad \sum_{i=1}^r w_i x_i^\ell = \frac{1}{2} \int_{-1}^1 x^\ell \, dx, \qquad \ell = 0, \ldots, z,$$

*if and only if the following two conditions are satisfied:*

(A) *The polynomial $V_r(x) = \prod_{i=0}^r (x - x_i)$ is orthogonal to all polynomials of degree $z - r$ with respect to the Lebesgue measure, that is,*

$$(3.3) \qquad \int_{-1}^1 V_r(x) x^\ell \, dx = 0, \qquad \ell = 0, \ldots, z - r.$$

(B) *The weights $w_j$ are given by*

$$(3.4) \qquad w_j = \frac{1}{2} \int_{-1}^1 \ell_j(x) \, dx,$$

*where*

$$\ell_i(x) = \prod_{j=1, j\neq i}^r \frac{x - x_j}{x_i - x_j}$$

*denotes the ith Lagrange interpolation polynomial with nodes $x_1, \ldots, x_r$.*

It follows from condition (3.3) that $z \leq 2r - 1$ (otherwise there does not exist a solution of this system). Moreover, there exists at least one set of points that satisfies condition (3.3), namely

$$(3.5) \qquad \{x \mid P_r(x) = 0\},$$

where $P_r$ denotes the $r$th Legendre polynomial orthogonal with respect to the Lebesgue measure on the interval $[-1, 1]$, because it is well known that this polynomial has $r$ distinct roots located in the interval $(-1, 1)$ (see [22]). Moreover, the corresponding weights defined by (3.4) are positive and define a Gaussian quadrature formula, that is, a quadrature formula of degree $2r-1$ with $r$ nodes (see, e.g., [6] or [8]). In the following discussion we will use this lemma with $z = 2d$. Then it follows that for any $r \in \{d+1, \ldots, 2d\}$, there exists at least one quadrature formula $\{x_i, w_i\}_{i=1}^r$ determined by the equations (3.3) and (3.4) that integrates polynomials of degree $2d$ exactly [in other words, the system of equations in (3.2) is satisfied with $z = 2d$] and the corresponding weights $w_i$ are positive. We consider a quadrature formula of this type. Define

$$(3.6) \qquad \theta_i = \arccos x_i, \qquad i = 1, \ldots, r,$$

and consider the design

$$(3.7) \qquad \mu = \begin{pmatrix} \theta_1 & \cdots & \theta_r \\ w_1 & \cdots & w_r \end{pmatrix}$$

on the interval $[0, \pi]$. Similarly, we define for any $t \in \mathbb{N}$ and any $\alpha \in (-\frac{t+1}{t}\pi, -\pi]$ a design $\nu = \nu(\alpha, t)$ on the interval $(-\pi, \pi]$ by

$$(3.8) \qquad \nu = \nu(\alpha, t) = \begin{pmatrix} \phi_1 & \cdots & \phi_t \\ \frac{1}{t} & \cdots & \frac{1}{t} \end{pmatrix},$$

where the points $\phi_j$ are given by

$$(3.9) \qquad \phi_j = \alpha + \frac{2\pi j}{t}, \qquad j = 1, \ldots, t.$$

The following result shows that designs of the form $\mu \otimes \nu$ are discrete $\Phi_p$- and $\Psi_{p,r}$-optimal designs for the spherical harmonic regression model (2.6).

THEOREM 3.2.   *Let $p \in [-\infty, 1)$, let $0 \le k_0 < \cdots < k_q \le d$ and denote by $K = K_{k_0, \ldots, k_q}$ the matrix defined by (2.12). For any $t \ge 2d+1$ and any $r \in \{d+1, \ldots, 2d\}$, the design $\mu \otimes \nu$ with factors given by (3.7) (corresponding to a quadrature formula of degree $2d$) and (3.8) is $\Phi_p$-optimal for estimating the coefficients $K^T c$, and is $\Psi_{p,r}$-optimal in the spherical harmonic regression model (2.6).*

PROOF.   Observing the proofs of Theorem 2.1 and Corollary 2.2, the assertion can be established by showing the identity

$$(3.10) \qquad M(\mu \otimes \nu) = I_{(d+1)^2}.$$

For this let

$$\psi(\phi) = (\psi_{-d}(\phi), \psi_{-d+1}(\phi), \dots, \psi_d(\phi))^T$$
$$= (\sqrt{2}\sin(d\phi), \dots, \sqrt{2}\sin\phi, 1, \sqrt{2}\cos\phi, \dots, \sqrt{2}\cos(d\phi))^T.$$

Then the regression functions in the spherical harmonic regression model of degree $d$ are given by

$$(3.11) \qquad \gamma_{ij} P_i^{|j|}(\cos\theta) \cdot \psi_j(\phi), \qquad j = -i, \dots, i,\, i = 0, \dots, d,$$

where $P_i^j$ is the Legendre function defined by (2.4) and the constants $\gamma_{ij}$ are given by

$$\gamma_{ij} = \sqrt{(2i+1)\frac{(i-|j|)!}{(i+|j|)!}}.$$

[Note that the different scaling of the cases $j = 0$ and $j \neq 0$ in (2.2) has been accommodated by introducing the factor $\sqrt{2}$ in the definition of the functions $\psi_j$.] Therefore, the identity (3.10) is equivalent to the system of equations

$$(3.12) \quad \gamma_{ij}\gamma_{k\ell} \int_{-\pi}^{\pi}\int_0^{\pi} P_i^{|j|}(\cos\theta)\psi_j(\phi)P_k^{|\ell|}(\cos\theta)\psi_\ell(\phi)\,d\mu(\theta)\,d\nu(\phi) = \delta_{ik}\delta_{j\ell}$$

$(j = -i, \dots, i; i = 0, \dots, d; \ell = -k, \dots, k; k = 0, \dots, d)$, where $\delta_{ik}$ denotes Kronecker's symbol. Observing Fubini's theorem, this system is satisfied if the equations

$$(3.13) \qquad \gamma_{ij}\gamma_{k\ell} \int_0^{\pi} P_i^j(\cos\theta)P_k^\ell(\cos\theta)\,d\mu(\theta) = \delta_{ik}\delta_{j\ell}$$

$(j = 0, \dots, i; i = 0, \dots, d; \ell = 0, \dots, k; k = 0, \dots, d)$ and

$$(3.14) \qquad \int_{-\pi}^{\pi} \psi_j(\phi)\psi_\ell(\phi)\,d\nu(\phi) = \delta_{j\ell}, \qquad j, \ell = -d, \dots, d,$$

can be established. It is well known (see Pukelsheim [16]) that the last identity is satisfied for measures of the form (3.8). To prove the remaining identity (3.13), recall that the measure $\mu$ corresponds to a quadrature formula that integrates polynomials of degree $2d$ exactly. Moreover, it follows from the representation (2.4) that for an even $m \in \{0, \dots, n\}$, the Legendre function $P_n^m(x)$ is a polynomial of degree $n$, while for any odd $m \in \{0, \dots, n\}$, the function $P_n^m(x)/\sqrt{1-x^2}$ is a polynomial of degree $n-1$. Observing the equation (3.14), we can restrict ourselves to the case $j = \ell$ for which the integrand $P_i^j(x)P_k^j(x)$ in (3.13) is always a polynomial of degree $i + k \leq 2d$, which contains the factor $(1 - x^2)$ if $j$ is odd. Consequently, we obtain from (3.2) (with $p = 2d$)

$$\int_0^\pi P_i^j(\cos\theta)P_k^j(\cos\theta)\,d\mu(\theta) = \sum_{\ell=1}^r w_\ell P_i^j(x_\ell)P_k^j(x_\ell)$$

$$= \frac{1}{2}\int_{-1}^1 P_i^j(x)P_k^j(x)\,dx = \frac{\delta_{ik}}{\gamma_{ij}\gamma_{kj}},$$

where the last property follows from the orthogonality of the Legendre functions (see [1]). This implies (3.12) [or equivalently (3.10)] and proves the theorem. □

REMARK 3.3.   It should be noted that the mapping $(\theta, \phi) \rightarrow (\sin\theta\cos\phi, \sin\theta\sin\phi, \cos\theta)$ from the rectangle $[0, \pi] \times (-\pi, \pi]$ onto the unit sphere $S_2$ maps all points of the form $(0, \phi)$ and $(\pi, \phi)$ with $\phi \in (-\pi, \pi]$ onto the points $(0, 0, 1)$ and $(0, 0, -1)$ on $S_2$, respectively. Moreover, it is easy to see that the vector $f_d(\theta, \phi)$ defined in (2.7) satisfies, for all $\phi \in (-\pi, \pi]$,

$$f_d(0, \phi) = (1, 0, \ldots, 0)^T \in \mathbb{R}^{(d+1)^2},$$

$$f_d(\pi, \phi) = (-1, 0, \ldots, 0)^T \in \mathbb{R}^{(d+1)^2}.$$

As a consequence, various points of the designs $\mu \otimes \nu$ constructed by Theorem 3.2 can be identified on the unit sphere if the support of the factor $\mu$ contains the point $0$ or $\pi$. To be precise, let $\bar{\mu}$ denote the measure obtained from $\mu$ by omitting the points $0$ and $\pi$, and define $\mu_0$ as the measure that puts masses $\mu(\{0\})$ and $\mu(\{\pi\})$ at the points $(0, 0)^T$ and $(\pi, 0)^T$, respectively. Then the measure

$$\mu_0 + \bar{\mu} \otimes \nu$$

has the same information matrix as the measure $\mu \otimes \nu$. Note that in the case $\mu(\{0\}) + \mu(\{\pi\}) = 0$ it follows that $\bar{\mu} = \mu$, because the points $0$ and $\pi$ are not support points of the design $\mu$.

EXAMPLE 3.4.   Consider the spherical harmonic regression model of degree $d = 1$ and the case $r = d + 1 = 2$. From Lemma 3.1 with $p = 2d = 2$ it follows that the polynomial $V_2(x) = (x - 1)(x + 1/3)$ satisfies

$$\int_{-1}^1 V_2(x)\,dx = 0.$$

The points $x_1 = -1/3$ and $x_2 = 1$ generate a quadrature formular of degree 2 on the interval $[-1, 1]$ with corresponding weights

$$w_1 = \int_{-1}^1 \frac{x-1}{(-4/3)}\frac{dx}{2} = \frac{3}{4}, \qquad w_2 = \int_{-1}^1 \frac{x+1/3}{4/3}\frac{dx}{2} = \frac{1}{4}.$$

According to Theorem 3.2 any design of the form

$$\mu \otimes \nu(\alpha, t) = \begin{pmatrix} 0 & \arccos(-\frac{1}{3}) \\ \frac{1}{4} & \frac{3}{4} \end{pmatrix} \otimes \nu(\alpha, t)$$

with $t \geq 3, \alpha \in (-\frac{t+1}{t}\pi, -\pi]$ is $\Phi_p$-optimal for estimating the parameters $K_{k_0, k_q}^T c$ ($q \leq 1$) and $\Psi_{p,r}$-optimal in the first-order spherical harmonic regression model (2.6). A typical example is given by the six-point design

$$\mu \otimes \nu = \begin{pmatrix} 0 & \arccos\left(-\frac{1}{3}\right) \\ \frac{1}{4} & \frac{3}{4} \end{pmatrix} \otimes \begin{pmatrix} -\frac{\pi}{3} & \frac{\pi}{3} & \pi \\ \frac{1}{3} & \frac{1}{3} & \frac{1}{3} \end{pmatrix},$$

and by Remark 3.3 the design with equal masses at the points $(\arccos(-\frac{1}{3}), -\frac{\pi}{3})$, $(\arccos(-\frac{1}{3}), \frac{\pi}{3})$, $(0,0)$ and $(\arccos(-\frac{1}{3}), \pi)$ has the same information matrix as the design $\mu \otimes \nu$, namely the identity matrix $I_4$.

REMARK 3.5.    Note that there are numerous possible ways to construct a discrete design with an information matrix equal to $I_{(d+1)^2}$ in the spherical harmonic regression model (2.6). According to Theorem 3.2, a quadrature formula with $r \in \{d+1, \ldots, 2d\}$ nodes $x_1, \ldots, x_r$ with positive weights is required, which integrates polynomials up to degree $2d$ exactly. By (3.7), this formula gives the factor $\mu$ of the optimal design $\mu \otimes \nu$, where the second factor is any design of the form (3.8) with $t \geq 2d+1$. By Lemma 3.1, the quadrature formula is determined by the equations

$$(3.15) \qquad \int_{-1}^{1} V_r(x) x^\ell \, dx = 0, \qquad \ell = 0, \ldots, 2d - r.$$

If $P_j(x)$ denotes the $j$th Legendre polynomial orthogonal with respect to the Lebesgue measure on the interval $[-1, 1]$, then the polynomial $V_r(x)$ can be represented as a linear combination of the Legendre polynomials $P_0(x), \ldots, P_r(x)$ and the orthogonality in (3.15) implies, for some constants $a_{2d-r+1}, \ldots, a_r$,

$$(3.16) \qquad V_r(x) = \sum_{j=2d-r+1}^{r} a_j P_j(x).$$

Note that the constants $a_{2d-r+1}, \ldots, a_r$ have to be chosen such that $V_r(x)$ has $r$ real roots in the interval $[-1, 1]$ and such that the corresponding weights defined by (3.4) are positive. This is in general a nontrivial problem. However, one can easily describe a class of quadrature formulas for which this property is satisfied. For this, let $P_j^{(\alpha, \beta)}(x)$ denote the $j$th Jacobi polynomial orthogonal with respect to the measure $(1-x)^\alpha (1+x)^\beta \, dx$ on the interval $[-1, 1]$ (see [22]). For any $r \geq d + 1$, it follows from these orthogonality properties that the identity (3.15) is satisfied for the polynomials

$$
\begin{array}{ll}
(3.17) & P_r^{(0,0)}(x), & (1-x)P_{r-1}^{(1,0)}(x), \\
& (1+x)P_{r-1}^{(0,1)}(x), & (1-x^2)P_{r-1}^{(1,1)}(x).
\end{array}
$$

Note that $P_r^{(0,0)}(x)$ is proportional to the Legendre polynomial $P_r(x)$ and consequently the representation of the form (3.16) is obvious, in this case choosing $a_{r-1} = \cdots = a_{2d-r+1} = 0$. Moreover, if $p$ is the degree of one of these polynomials, it follows from classical results on orthogonal polynomials that each of the polynomials in (3.17) has precisely $p$ roots in the interval $[-1, 1]$. It is shown in [2] that the weights corresponding, by formula (3.4), to these roots are positive (see also [8]), and as a consequence we obtain a quadrature formula of degree $2d$ on the interval $[-1, 1]$, say $\{x_j, w_j\}_{j=1,\ldots,p}$. The corresponding design $\mu$ obtained from (3.6) and (3.7) gives, in combination with any design $\nu$ of the form (3.8), a $\Phi_p$-optimal design $\mu \otimes \nu$ for estimating the parameters $K_{k_1,\ldots,k_q}^T c$ in the spherical harmonic regression model (2.6). This design is also $\Psi_{p,r}$-optimal by Corollary 2.2. We finally note that the zeros of the polynomials in (3.17) are the support points of $D$-optimal designs in heteroscedastic polynomial regression models (see [20]).

Although there are numerous designs on the rectangle $[0, \pi] \times (-\pi, \pi]$ with information matrix in the spherical harmonic regression model (2.6) given by $I_{(d+1)^2}$, the support points of the factor corresponding to the polar angle $\theta$ have to cover a sufficiently large range of the interval $[0, \pi]$ in order to obtain a $\Phi_p$-optimal design in the sense of Theorem 3.2. This is the statement of the final theorem of this section, which will be proven in the Appendix.

THEOREM 3.6. *Let $x_1^*$ denote the smallest zero of the Legendre polynomial $P_{d+1}(x)$ and define $z^* = \arccos|x_1^*|$. If $z > z^*$, then there exists no design on the rectangle $[z, \pi - z] \times (-\pi, \pi]$ with information matrix $I_{(d+1)^2}$ in the spherical harmonic regression model (2.6).*

## 4. Further discussion.

4.1. *The second-order spherical harmonic model.* Consider the case $d = 2$ corresponding to the second-order spherical harmonic model with nine regression functions. We calculate the designs corresponding to the four cases in (3.17) for the choice $r = d + 1 = 3$. From $P_3(x) = x(5x^2 - 3)/2$ we obtain the support points of the probability measure $\mu$ corresponding to the polar angle $\theta$ as

$$
(4.1) \quad \theta_1 = \arccos\sqrt{\frac{3}{5}}, \qquad \theta_2 = \arccos 0 = \frac{\pi}{2}, \qquad \theta_3 = \arccos\left(-\sqrt{\frac{3}{5}}\right),
$$

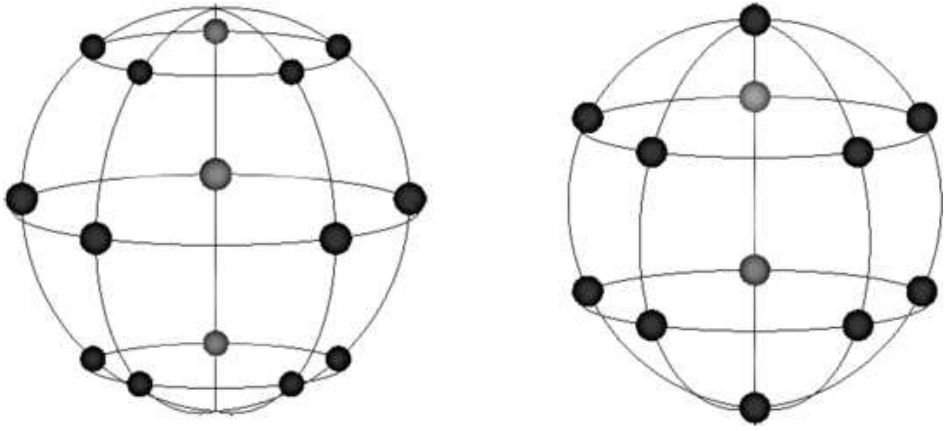

Fig. 3.    *Projections of the support points of the optimal design in the second-order spherical harmonic regression model* (2.6). *The left panel corresponds to the optimal design defined by* (4.1), *while the right panel represents the design* (4.3).

and the corresponding weights are given by

$$(4.2) \qquad w_2 = \frac{1}{2}\int_{-1}^{1} \frac{x^2 - 3/5}{(-3/5)}\,dx = \frac{4}{9}, \qquad w_1 = w_3 = \frac{1-w_2}{2} = \frac{5}{18}.$$

Similarly, if the polynomial $(x^2-1)P_2^{(1,1)}(x) = \frac{3}{4}(x^2-1)(5x^2-1)$ is used for the construction of a quadrature formula, we obtain

$$(4.3) \qquad \begin{aligned} &\theta_1 = \arccos 1 = 0, &&\theta_2 = \arccos\sqrt{\frac{1}{5}}, \\[2mm] &\theta_3 = \arccos\left(-\sqrt{\frac{1}{5}}\right), &&\theta_4 = \arccos(-1) = \pi, \end{aligned}$$

and the weights are obtained from the representation (3.4) and given by

$$(4.4) \qquad w_1 = w_4 = \frac{1}{12}, \qquad w_2 = w_3 = \frac{5}{12}.$$

The measure $\nu$ corresponding to the azimuthal angle $\phi$ is given by (3.8), where $t \geq 5$. The projections of the support points of the two measures onto the unit sphere are depicted by the left and right panels of Figure 3, where for the second component a design with $t = 5$ support points and $\alpha = -\pi$ is used. It should be noted that it follows from Remark 3.3 that the design $\mu \otimes \nu$ obtained from (4.3) and (4.4) for the factor $\mu$ and $\nu = \nu(\alpha, 5)$ corresponds to uniform design on the sphere with 12 support points. A general construction of uniform designs will be discussed in the following paragraph. We finally illustrate the two other non-symmetric cases in (3.17). For the supporting

polynomial we obtain

$$(1-x)P_2^{(1,0)}(x) = \frac{1-x}{2}(5x^2 + 2x - 1),$$

which gives, for the support points of the factor $\mu$ of the product design $\mu \otimes \nu(\alpha, t)$,

$$
\begin{aligned}
(4.5) \qquad
\theta_1 &= \arccos 1 = 0, \\
\theta_2 &= \arccos\left(\frac{-1+\sqrt{6}}{5}\right), \\
\theta_3 &= \arccos\left(\frac{-1-\sqrt{6}}{5}\right),
\end{aligned}
$$

with corresponding weights

$$w_1 = \tfrac{1}{9}, \qquad w_2 = \tfrac{1}{36}(16+\sqrt{6}), \qquad w_3 = \tfrac{1}{36}(16-\sqrt{6}).$$

The fourth case of a factor $\mu$ of a design $\mu \otimes \nu(\alpha, t)$ yielding $I_9$ as information matrix in the second-order spherical harmonic model (2.6) is obtained by symmetry, that is,

$$
\begin{aligned}
(4.6) \qquad
\theta_1 &= \arccos\left(\frac{1+\sqrt{6}}{5}\right), \\
\theta_2 &= \arccos\left(\frac{1-\sqrt{6}}{5}\right), \qquad \theta_3 = \arccos(1) = \pi, \\
w_1 &= \tfrac{1}{36}(16-\sqrt{6}), \qquad w_2 = \tfrac{1}{36}(16+\sqrt{6}), \qquad w_3 = \tfrac{1}{9}.
\end{aligned}
$$

The projections of the support points onto the unit sphere of the corresponding product designs $\mu \times \nu(\alpha, t)$ with $\alpha = -\pi$ and $t = 5$ are depicted by the left and right panels of Figure 4. Note that the two cases are related by a reflection of the support points at the equator.

4.2. *Optimal designs with equal weights.* Note that the designs provided by Theorem 3.2 are in general not uniform designs, which would have equal weights at their support points. Because designs with this structure are particularly attractive from a practical point of view, we will briefly discuss the possibility of their construction in this section. Note that for the determination of an optimal design in the spherical harmonic regression model (2.6) with equal weights at its supports by the procedure introduced in Section 3, it is necessary to find a quadrature formula that has equal weights at its support points and integrates polynomials of degree $2d$ exactly. This problem has a long history in mathematics (see [6]).

It follows from Lemma 3.1 that such a formula must have at least $d+1$ nodes. Moreover, quadrature formulas with the minimal number of $n$ nodes and equal weights which integrate polynomials of degree $n$ exactly only exist in the cases $n = 0, 1, \ldots, 7, 9$ (see [6], page 58). In the present context these formulas correspond to uniform designs for the polar angle $\theta$ in the case $d = 1, 2, 3, 4$, which give, in combination with a design of the form (3.8), a $\Phi_p$- or $\Psi_{p,r}$-optimal design in the spherical harmonic regression model (2.6) with equal weights at its support points. The corresponding nodes of the quadrature formula required for the factor $\mu$ are depicted in the first four rows of Table 1. If $d \geq 5$, quadrature formulas with equal weights on $2d$ points integrating polynomials of degree $2d$ exactly do not exist and more nodes are required for the construction of such formulas. We determined formulas of this type numerically for $d = 5, 6, 7$ and depicted them in the last three rows of Table 1. Note that these formulas use the origin and that the number of nodes increases rapidly. For example, if $d = 7$, we only found a quadrature formula with 23 support points and equal weights which integrates polynomials of degree 14 exactly.

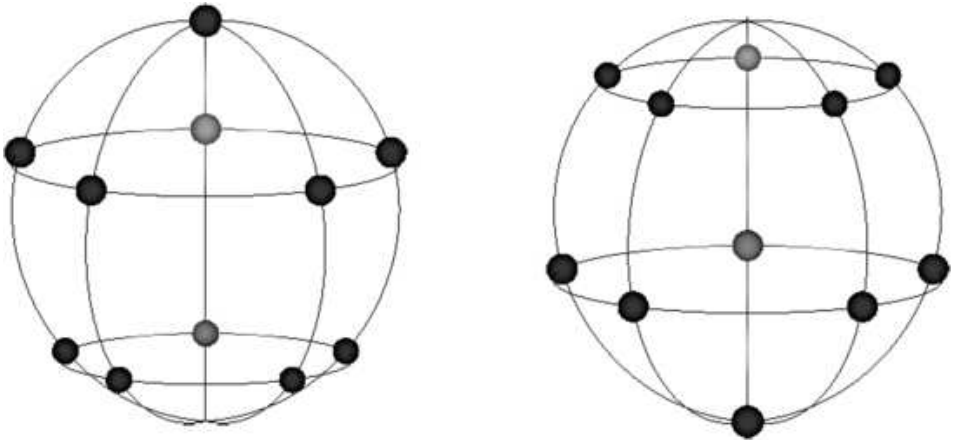

FIG. 4. *Projections of the support points of the optimal design in the second-order spherical harmonic regression model (2.6). The left panel corresponds to the optimal design defined by (4.5), while the right panel represents the design (4.6).*

4.3. *Some efficiency considerations.* It might be of interest to compare some of the commonly used designs from the literature on the analysis of 3D shapes with the $\Psi_{p,r}$- and $\Phi_p$-optimal designs obtained in the present paper. As an example, we consider two uniform designs of the form $\mu \otimes \nu$, where the components are given by

$$(4.7) \qquad \mu = \begin{pmatrix} \dfrac{i}{n_1 + 1} \\ \dfrac{1}{n_1} \end{pmatrix}_{i=1,\ldots,n_1}, \qquad \nu = \begin{pmatrix} \dfrac{2i\pi}{n_2} - \pi \\ \dfrac{1}{n_2} \end{pmatrix}_{i=1,\ldots,n_2},$$

TABLE 1
*Quadrature formulas with equal weights at their nodes*[*]

| $d$ | $x_0$ | $\pm x_1$ | $\pm x_2$ | $\pm x_3$ | $\pm x_4$ | $\pm x_5$ | $\pm x_6$ | $\pm x_7$ | $\pm x_8$ | $\pm x_9$ | $\pm x_{10}$ | $\pm x_{11}$ |
|---|---|---|---|---|---|---|---|---|---|---|---|---|
| 1 | — | 0.577 | — | — | — | — | — | — | — | — | — | — |
| 2 | — | 0.188 | 0.795 | — | — | — | — | — | — | — | — | — |
| 3 | — | 0.267 | 0.423 | 0.866 | — | — | — | — | — | — | — | — |
| 4 | 0 | 0.168 | 0.529 | 0.601 | 0.912 | — | — | — | — | — | — | — |
| 5 | 0 | 0.223 | 0.247 | 0.443 | 0.671 | 0.724 | 0.939 | — | — | — | — | — |
| 6 | 0 | 0.008 | 0.282 | 0.358 | 0.458 | 0.566 | 0.760 | 0.778 | 0.954 | — | — | — |
| 7 | 0 | 0.174 | 0.177 | 0.186 | 0.328 | 0.502 | 0.533 | 0.542 | 0.712 | 0.797 | 0.852 | 0.965 |

[*]These formulas correspond to a uniform design $\mu_u$ of the form (3.7) for the polar angle $\theta$. The resulting design $\mu_u \otimes \nu(\alpha, t)$ with $\nu(\alpha, t)$ defined by (3.8), $t \geq 2d + 1$, is a $\Phi_p$- or $\Psi_{p,r}$-optimal design in the spherical harmonic regression model (2.6) of degree $d$ with equal weights at its support points.

and the designs are given by

$$(4.8) \quad \mu = \begin{pmatrix} \arccos\left(1 - \dfrac{2i}{n_1 + 1}\right) \\ \dfrac{1}{n_1} \end{pmatrix}_{i=1,\ldots,n_1}, \qquad \nu = \begin{pmatrix} \dfrac{2i\pi}{n_2} - \pi \\ \dfrac{1}{n_2} \end{pmatrix}_{i=1,\ldots,n_2},$$

respectively. Note that the design $\mu \otimes \nu$ defined by (4.7) corresponds to a uniform distribution on the grid in the rectangle $[0, \pi] \times [-\pi, \pi]$, while a design $\mu \otimes \nu$ of the form (4.8) yields a design on the sphere that takes observations on several circles with equal distance on the $z$ axis. This design was used by Ding, Nesumi, Takano and Ukai [5] for a principal component anal-

TABLE 2
*D, E, A and $\Psi_{-1,r}$ efficiencies ($r = 2, 3$) of the uniform designs (4.7) and (4.8) in first- and second-order spherical harmonic models*

| | | Design (4.7) | | | | | Design (4.8) | | | | |
|---|---|---|---|---|---|---|---|---|---|---|---|
| $d$ | $n_1$ | $\mathrm{eff}_D$ | $\mathrm{eff}_E$ | $\mathrm{eff}_A$ | $\mathrm{eff}_{\Psi_{-1,2}}$ | $\mathrm{eff}_{\Psi_{-1,3}}$ | $\mathrm{eff}_D$ | $\mathrm{eff}_E$ | $\mathrm{eff}_A$ | $\mathrm{eff}_{\Psi_{-1,2}}$ | $\mathrm{eff}_{\Psi_{-1,3}}$ |
| 1 | 3 | 1.000 | 1.000 | 1.000 | 1.000 | 1.000 | 0.940 | 0.500 | 0.870 | 0.667 | 0.789 |
| | 4 | 0.997 | 0.938 | 0.994 | 0.938 | 0.957 | 0.964 | 0.600 | 0.923 | 0.750 | 0.857 |
| | 5 | 0.993 | 0.900 | 0.986 | 0.900 | 0.931 | 0.976 | 0.667 | 0.949 | 0.800 | 0.894 |
| | 6 | 0.989 | 0.875 | 0.979 | 0.875 | 0.913 | 0.983 | 0.714 | 0.964 | 0.833 | 0.916 |
| | 7 | 0.986 | 0.857 | 0.973 | 0.857 | 0.900 | 0.987 | 0.750 | 0.973 | 0.857 | 0.931 |
| 2 | 4 | 0.991 | 0.801 | 0.982 | 0.838 | 0.851 | 0.902 | 0.229 | 0.745 | 0.331 | 0.427 |
| | 5 | 0.987 | 0.805 | 0.974 | 0.824 | 0.831 | 0.935 | 0.323 | 0.838 | 0.435 | 0.539 |
| | 6 | 0.981 | 0.794 | 0.964 | 0.807 | 0.811 | 0.954 | 0.399 | 0.888 | 0.512 | 0.617 |
| | 7 | 0.976 | 0.782 | 0.955 | 0.793 | 0.796 | 0.965 | 0.461 | 0.918 | 0.571 | 0.674 |
| | 8 | 0.972 | 0.772 | 0.947 | 0.782 | 0.785 | 0.973 | 0.511 | 0.937 | 0.617 | 0.716 |

Table 3

*D, E, A and $\Psi_{-1,r}$ efficiencies ($r = 2, 3$) of the designs (4.7) and (4.8) in spherical harmonic models of degree 3 and 4*

| | | Design (4.7) | | | | | Design (4.8) | | | | |
|---|---|---|---|---|---|---|---|---|---|---|---|
| $d$ | $n_1$ | eff$_D$ | eff$_E$ | eff$_A$ | eff$_{\Psi_{-1,2}}$ | eff$_{\Psi_{-1,3}}$ | eff$_D$ | eff$_E$ | eff$_A$ | eff$_{\Psi_{-1,2}}$ | eff$_{\Psi_{-1,3}}$ |
| 3 | 5 | 0.980 | 0.799 | 0.961 | 0.802 | 0.803 | 0.874 | 0.094 | 0.600 | 0.146 | 0.199 |
| | 6 | 0.975 | 0.784 | 0.953 | 0.784 | 0.787 | 0.911 | 0.155 | 0.733 | 0.223 | 0.296 |
| | 7 | 0.970 | 0.768 | 0.944 | 0.768 | 0.772 | 0.934 | 0.214 | 0.811 | 0.292 | 0.377 |
| | 8 | 0.965 | 0.756 | 0.936 | 0.756 | 0.760 | 0.948 | 0.269 | 0.859 | 0.352 | 0.444 |
| | 9 | 0.961 | 0.747 | 0.929 | 0.747 | 0.751 | 0.959 | 0.318 | 0.891 | 0.404 | 0.500 |
| 4 | 6 | 0.969 | 0.739 | 0.942 | 0.755 | 0.761 | 0.851 | 0.035 | 0.434 | 0.057 | 0.081 |
| | 7 | 0.965 | 0.747 | 0.936 | 0.751 | 0.753 | 0.890 | 0.067 | 0.600 | 0.102 | 0.142 |
| | 8 | 0.961 | 0.739 | 0.929 | 0.742 | 0.743 | 0.916 | 0.103 | 0.709 | 0.149 | 0.203 |
| | 9 | 0.956 | 0.731 | 0.922 | 0.733 | 0.734 | 0.933 | 0.142 | 0.781 | 0.196 | 0.261 |
| | 10 | 0.952 | 0.724 | 0.915 | 0.726 | 0.727 | 0.945 | 0.180 | 0.830 | 0.240 | 0.314 |

ysis of the variance–covariance matrix of the least squares estimator (1.3) in the spherical harmonic regression model of degree 7. In Tables 2 and 3 we consider the problem of designing an experiment to estimate the full parameter vector $c$ or the principal component analysis in the spherical harmonic regression model (2.6). We show the $D$, $E$ and $A$ efficiencies of these designs for various values of $n_1$ and $n_2$ in the spherical harmonic models of degree 1, 2, 3 and 4. The tables also contain the $\Psi_{-1,r}$-efficiencies, which are defined by

$$(4.9) \qquad \text{eff}_{\Psi,p,r}(\xi) = \frac{\Psi_{p,r}(\xi)}{\sup_\eta \Psi_{p,r}(\eta)}.$$

For the sake of brevity we choose the design for the azimuthal angle $\phi$ as the design defined by (3.8)—other uniform designs for this component will yield substantially lower efficiencies and are therefore not depicted. As a consequence, the efficiencies of the design $\mu \otimes \nu$ do not depend on $n_2$ (provided that $n_2 \geq 2d + 1$, which will be assumed throughout this section).

For the first-order spherical harmonic regression model, we observe very good $D$ and $A$ efficiencies of the designs (4.7) and (4.8). However, the $E$ efficiencies and the $\Psi_{p,r}$ efficiencies of these designs are substantially smaller, in particular those efficiencies obtained for the design (4.8) with moderate values of $n_1$. For spherical harmonic models of larger degree both designs will still yield high $D$ efficiencies, the designs defined by (4.7) yield reasonable $A$ efficiencies, but the $E$ and $\Psi_{-1,r}$ efficiencies are substantially smaller. Moreover, the $A$ and $E$ efficiencies of the design (4.8) are very low. The $\Psi_{-1,r}$ efficiencies of this design are slightly larger but still not satisfactory. It is also interesting to note that the efficiencies of the design (4.7) are decreasing with the number $n_1$ while they are increasing for the design (4.8).

From these calculations and additional results, which are not depicted for the sake of brevity, we observe that designs with a uniform distribution on a grid in the rectangle $[0, \pi] \times (-\pi, \pi]$ should be used only if maximization of the $D$ criterion is the preliminary goal of the design of the experiment. Whenever principal component analysis is the main goal of the experiment or precise estimates of the parameters themselves are required, some more care is necessary in the design of an experiment for analyzing 3D shapes. In this case, the loss of efficiency when using uniform designs on a grid in the rectangle $[0, \pi] \times [-\pi, \pi]$ or a uniform design taking observations on several circles with equal distance on the $z$ axis can be substantial. In most cases there exist substantially more efficient designs for the analysis in a spherical harmonic regression model. The advantages of the optimal designs derived in the present paper will also be illustrated in the following example.

4.4. *A concluding example.* To demonstrate the benefits of our designs, we finally reanalyze the design used by Ding, Nesumi, Takano and Ukai [5] for shape analysis of Citrus species. These authors used observations at 360 points measured in 10 circles using the equal height sampling method for the $z$ axis. By choosing $d = 7$ the data for the surface shape were expanded into the first 64 terms of spherical harmonic functions. The information contained in the spherical harmonic coefficients was summarized by a principal component analysis using the first seven principal components. Note that the design of Ding, Nesumi, Takano and Ukai [5] corresponds to a design of the form (4.8) with $n_1 = 10$ and $n_2 = 36$. In the first row of Table 4 we show the efficiencies of this design with respect to the optimal designs obtained in this paper. The first factor of an optimal (uniform) design can be obtained from the quadrature formula corresponding to the spherical harmonic regression model of degree 7 in Table 1. We observe a reasonable efficiency only for the $D$-criterion. For all other criteria the design used by Ding, Nesumi, Takano and Ukai [5] is very inefficient. The optimal designs proposed in this paper (or appropriate approximations) will yield substantially smaller variances of the least squares estimates for linear combinations of the parameters.

We finally note that the optimal designs obtained in this paper are approximate and do not have masses that are multiples of $1/360$. However, a very efficient design for the inference in the spherical harmonic regression model of degree 7 using 360 different points can easily be obtained as follows. The quadrature formula obtained from Table 1 has equal masses at 23 nodes and yields a design of the form (3.7), say $\{\theta_i, 1/23\}_{i=1}^{23}$, where $\theta_i = \arccos x_i$ and $x_1, \ldots, x_{23}$ are the points listed in Table 1 for the case $d = 7$. We propose to combine this design with two designs of the form (3.8) to obtain an efficient exact design. More precisely, we propose to use the *D, A, E and efficiencies* $\mathrm{eff}_{\Psi_{-1,r}}$ $(r = 1, \dots, 10)$ *in the spherical harmonic regression model of degree 7 used by Ding, Nesumi, Takano and Ukai* [5][*]

| | | | | $\mathrm{eff}_{\Psi_{-1,r}}$ | | | | | | | |
|---|---|---|---|---|---|---|---|---|---|---|---|
| **1** | **2** | **3** | **4** | **5** | **6** | **7** | **8** | **9** | **10** | $\mathrm{eff}_A$ | $\mathrm{eff}_D$ |
| 0.003 | 0.006 | 0.008 | 0.011 | 0.013 | 0.016 | 0.019 | 0.021 | 0.024 | 0.026 | 0.149 | 0.840 |
| 0.958 | 0.958 | 0.958 | 0.958 | 0.958 | 0.958 | 0.958 | 0.958 | 0.958 | 0.958 | 0.987 | 0.992 |

[*]The first row is the design given by (4.8) with $n_1 = 10$ and $n_2 = 36$. The second row gives the efficiencies of the design defined by (4.10), which takes one observation at each point of the set $\mathcal{U}$.

uniform distribution $\xi_U^*$ at the points

$$(4.10) \quad \begin{aligned} \mathcal{U} = &\left\{ \left( \theta_i, \frac{\pi(2j-15)}{15} \right) \Big| i = 1, \dots, 8; j = 1, \dots, 15 \right\} \\ &\cup \left\{ \left( \theta_i, \frac{\pi(2j-16)}{16} \right) \Big| i = 9, \dots, 23; j = 1, \dots, 16 \right\} \end{aligned}$$

as the design for the spherical harmonic regression model of degree 7. The implementation of this design is illustrated in Figure 5. The efficiencies of the exact design $\xi_U^*$ are depicted in the second row of Table 4 and we observe that this design, which advises the experimenter to take one observation at each point of $\mathcal{U}$, is highly efficient with respect to all optimality criteria under consideration.

**5. Conclusions and some directions for future research.** In this paper we tried to provide a starting point for studying design problems for 3D shape analysis. There are many issues in this area which can be addressed by the choice of an experimental design and we concentrated on spherical harmonic descriptors, with special emphasis on classification. This direction was motivated by recent work of Ding, Nesumi, Takano and Ukai [5], who used the coefficients in a spherical harmonic expansion to classify different fruit shapes by principal component analysis. Other authors compare different shapes by simply computing the Euclidean distance between their corresponding spherical harmonic descriptors (see [7]). If classification is the main object in 3D shape analysis, precise estimation of the coefficients in the spherical harmonic expansion is of particular importance and should be addressed by the choice of an experimental design, whenever this is possible. In the present paper we have constructed optimal designs for this purpose by considering Kiefer's $\Phi_p$ criteria and a new criterion which is directly related to principal component analysis. It is demonstrated that the new designs are substantially more efficient for parameter estimation and principal component analysis than the commonly used designs by reanalyzing a typical data example from the literature.

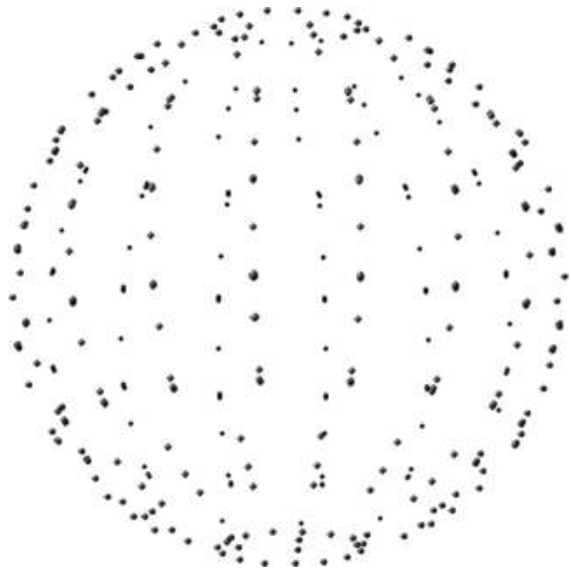

Fig. 5.   *Implementation of the design* (4.10) *on the unit sphere.*

As pointed out by a referee, there are several other design issues concerning image analysis that are not addressed in this work, but are important topics for future research in this direction. In the following list we briefly mention the—in our opinion—most important directions in image analysis, where the generalization of our approach could be profitable.

1. The main motivation of our work stems from the fact that spherical harmonic descriptors are used for classification by principal component analysis. There are many other statistical tools for this purpose, for example, independent component analysis (see the recent monograph by Hyvärinen, Karhunen and Oja [10]). Therefore, it is important to determine efficient designs for alternative classification rules and to compare these with the designs derived in this paper.

2. In addition to the problem of classification, another important goal of 3D shape analysis is the reconstruction of shapes. In such cases the number of basis functions is usually very large and an optimality criterion should focus on the closeness between the original shape and the reconstructed shape. A suitable optimality criterion for this purpose is to minimize the expected mean squared error with respect to the choice of an experimental design. This problem is closely related to the question of model uncertainty. Since the pioneering work of Box and Draper [3], several authors have addressed the problem of incorporating the bias in the optimality criterion (see, e.g., Wiens [24] and the references therein). We note, however, that by Theorem 2.1 the uniform distribution on the sphere is $\Phi_p$-optimal in the model (2.6) for any order $d \in \mathbb{N}$. Therefore, from a theoretical viewpoint, the question of model uncertainty does not exist for the uniform distribution on the sphere (it is the optimal design for any $d \in \mathbb{N}$). However, these problems arise if a design has to be implemented for a fixed sample size. If the main goal of the analysis is the reconstruction of shapes, we recommend choosing $d$ as large as possible. In Section 4.4 we have indicated how such a design can be found in the case $d = 7$, which was given in a concrete application of spherical harmonic descriptors in biology. In principle this method can be adapted to larger values of $d$. However, the numerical difficulties of finding quadrature formulas of degree $2d$ with equal weights increase substantially with the value $d$ and suitable software has to be developed for this purpose. We also note that the sample size $n$ always provides an upper bound for the degree of the model (2.6). Therefore—in principle—the problem of model uncertainty is caused by a limited sample size, not by our approach of constructing optimal designs.

3. It is remarkable that the approach developed in this paper does not require any distributional assumptions, because it is based on the least squares method and the covariance matrix of the corresponding estimator. However, the assumption of uncorrelated observations is important for the calculation of this matrix. Observations from image analysis are usually taken from one object and therefore this assumption has either to be checked or to be taken into account by the construction of the design. Note that many authors working in image analysis do not consider an error in distribution. In our model (2.6), the error refers to a nonsystematic measurement error, for which the assumption of no correlation may be justified. In general, the construction of optimal designs for models with correlated observations is a very hard problem even in less sophisticated models than considered in this paper (see, e.g., [18]). The incorporation of correlation structures in the construction of optimal designs is probably one of the most challenging problems in future research on optimal designs for image analysis.

4. It seems to be desirable to extend the methods derived in this paper to other problems in image analysis, in particular to the problem of 2D image reconstruction. We note that the design problem for 3D shape analysis is substantially simpler than design problems for general image analysis. The reason for this is the specific structure of the system of spherical harmonic functions used in (1.1), which is not available for 2D image analysis. A reasonable system for the two-dimensional case is the Zernike polynomials (see [14]), and research to construct good designs in this case is still ongoing. Initial results in this direction indicate that there are substantial differences between the two- and the three-dimensional cases.

## APPENDIX: TECHNICAL DETAILS

PROOF OF LEMMA 3.1.    Assume that conditions (A) and (B) of Lemma 3.1 are satisfied and let $Q(x)$ denote an arbitrary polynomial of degree $z$. By Bezout's theorem the polynomial $Q$ can be represented in the form

$$(A.1) \qquad Q(x) = P(x)V_r(x) + R(x),$$

where $V_r(x) = \prod_{j=1}^{r}(x - x_j)$, the polynomial $P(x)$ is of degree $z - r$ and the degree of $R(x)$ is less than $r$. Because the degree of $R$ is at most $r - 1$, it can be represented as

$$(A.2) \qquad R(x) = \sum_{j=1}^{r} \ell_j(x) R(x_j),$$

and we obtain from the conditions (A) and (B) of the lemma that

$$\frac{1}{2} \int_{-1}^{1} Q(x)\,dx = \frac{1}{2} \int_{-1}^{1} R(x)\,dx$$

$$= \frac{1}{2} \sum_{j=1}^{r} R(x_j) \int_{-1}^{1} \ell_j(x)\,dx$$

$$= \sum_{j=1}^{r} R(x_j) w_j = \sum_{j=1}^{r} Q(x_j) w_j.$$

Using the functions $Q(x) = x^\ell$ $(\ell = 0, \ldots, z - r)$ yields the identities in (3.2).

For a proof of the converse, we assume that (3.2) is valid and obtain for $\ell = 0, \ldots, z - r$,

$$\int_{-1}^{1} V_r(x) x^\ell\,dx = \sum_{j=1}^{r} V_r(x_j) x_j^\ell w_j = 0,$$

which gives condition (A). On the other hand, condition (B) follows from the identity

$$\frac{1}{2} \int_{-1}^{1} \ell_j(x)\,dx = \sum_{i=1}^{r} w_j \ell_j(x_i) = w_j,$$

observing the property $\ell_j(x_i) = \delta_{ij}$ for the Lagrange interpolation polynomials.    $\square$

PROOF OF THEOREM 3.6. Let $\xi$ denote an arbitrary discrete design on the rectangle $[0, \pi] \times (-\pi, \pi]$ and denote by $\theta_1, \ldots, \theta_N$ the distinct first coordinates of the corresponding support points of $\xi$. Obviously $\xi$ can be represented as

$$(A.3) \qquad \xi = \sum_{i=1}^{N} \xi_i,$$

where the designs $\xi_i$ are defined by

$$\xi_i = \begin{pmatrix} (\theta_i, \phi_{i1}) & \cdots & (\theta_i, \phi_{iN_i}) \\ w_{i1} & \cdots & w_{iN_i} \end{pmatrix}$$

(with $N_i \in \mathbb{N}$) and represent the part of $\xi$ corresponding to the support points with first coordinate equal to $\theta_i$ $(i = 1, \ldots, N)$. Define $x_i = \cos \phi_i$ and $w_i = \sum_{j=1}^{N_i} w_{ij}$ $(i = 1, \ldots, N)$, and consider the design

$$(A.4) \qquad \eta = \eta_\xi = \begin{pmatrix} x_1 & \cdots & x_N \\ w_1 & \cdots & w_N \end{pmatrix}.$$

If the design $\xi$ satisfies the condition $M(\xi) = I_{(d+1)^2}$, it follows for the submatrix corresponding to the $(d+1)$ regression functions $Y_\ell^0(\theta, \phi) = \sqrt{2\ell + 1} P_\ell^0(\cos \theta)$ $(\ell = 0, \ldots, d)$ that

$$
\begin{aligned}
I_{d+1} &= \left( \int_0^\pi \int_{-\pi}^\pi Y_\ell^0(\theta, \phi) Y_k^0(\theta, \phi) \, d\xi(\theta, \phi) \right)_{\ell, k = 0}^{d} \\
(A.5) \qquad &= \left( \sum_{i=1}^N \sum_{j=1}^{N_i} w_{ij} \sqrt{2\ell + 1} \sqrt{2k + 1} P_\ell^0(x_i) P_k^0(x_i) \right)_{\ell, k = 0}^{d} \\
&= \left( \sqrt{2\ell + 1} \sqrt{2k + 1} \sum_{i=1}^N w_i P_\ell(x_i) P_k(x_i) \right)_{\ell, k = 0}^{d}.
\end{aligned}
$$

[Note that $P_\ell^0(x) = P_\ell(x)$ is the $\ell$th Legendre polynomial.] On the other hand, the orthogonality relation for the Legendre polynomials (see [22]) yields

$$(A.6) \qquad \frac{1}{2} \int_{-1}^1 P_\ell(x) P_k(x) \, dx = \frac{\delta_{\ell k}}{2\ell + 1}, \qquad \ell, k = 0, \ldots, d,$$

and it follows from (A.5) that the quadrature formula defined by the design $\eta$ in (A.4) integrates polynomials of degree $2d$ exactly. The assertion of Theorem 3.6 now follows from the following auxiliary lemma. □

LEMMA A.1. *For any $n \in \mathbb{N}$, let $\Xi$ denote the set of all probability measures on the interval $[-1, 1]$ of the form (A.4) which integrate polynomials* *of degree $n$ exactly, that is,*

$$
(A.7) \qquad \tfrac{1}{2} \int_{-1}^{1} x^{\ell}\, dx = \sum_{j=1}^{n} w_j x_j^{\ell}, \qquad \ell = 0, \ldots, n.
$$

*Then*

$$
\min_{\eta \in \Xi} \min\{v \in \mathbb{R}^{+} \mid \operatorname{supp} \eta \subset [-v, v]\} = |x_1^{*}|,
$$

*where $x_1^{*}$ is the smallest zero of the $r$th Legendre polynomial $P_r(x)$ and $r = \lfloor n/2 \rfloor + 1$. Moreover, the minimum is attained by the measure*

$$
(A.8) \qquad \eta^{*} = \begin{pmatrix} x_1^{*} & \cdots & x_r^{*} \\ w_1^{*} & \cdots & w_r^{*} \end{pmatrix},
$$

*where the points $x_1^{*}, \ldots, x_r^{*}$ are the zeros of the Legendre polynomial $P_r(x)$, $w_i = \tfrac{1}{2} \int_{-1}^{1} \ell_i(x)\, dx$ $(i = 1, \ldots, r)$ and $\ell_i(x)$ is the Lagrange interpolation polynomial with nodes $x_1^{*}, \ldots, x_r^{*}$.*

PROOF. For a design $\eta$ on the interval $[-1, 1]$, let

$$
(A.9) \qquad u(\eta) := \min\{v \in \mathbb{R}^{+} \mid \operatorname{supp} \eta \subset [-v, v]\}
$$

denote the half of the minimal length of intervals that contain the support of $\eta$ and define for each $N \in \mathbb{N}$ the set $\Xi_N \subset \Xi$ as the set of all probability measures with $N$ support points satisfying (A.7) (note that $\Xi = \bigcup_{N \in \mathbb{N}} \Xi_N$). Because for a design $\eta$ of the form (A.4) the equations in (A.7) can be written as

$$
(A.10) \qquad \tfrac{1}{2} \int_{-1}^{1} \begin{pmatrix} 1 \\ \vdots \\ x^n \end{pmatrix} dx = \sum_{j=1}^{N} w_j \begin{pmatrix} 1 \\ \vdots \\ x_j^n \end{pmatrix},
$$

it follows from Carathéodory's theorem (see [19]) that there exists a measure, say $\tilde{\eta}$, with at most $n+1$ support points $x_{j_1}, \ldots, x_{j_{n+1}}$ $(1 \le j_1 < j_2 < \cdots < j_{n+1} \le N)$ such that (A.7) is satisfied. Consequently, we obtain

$$
(A.11) \qquad \inf_{\eta \in \Xi} u(\eta) = \inf_{\eta \in \Xi_{n+1}} u(\eta).
$$

If

$$
\eta = \begin{pmatrix} x_1 & \cdots & x_{n+1} \\ w_1 & \cdots & w_{n+1} \end{pmatrix} \in \Xi_{n+1}
$$

with $-1 \le x_1 < \cdots < x_{n+1} \le 1$ is any probability measure on the interval $[-1, 1]$ that satisfies (A.7), then by Lemma 3.1 the weights can be represented as

$$
w_i = \tfrac{1}{2} \int_{-1}^{1} \ell_i(x)\, dx,
$$

where $\ell_i(x)$ is the $i$th Lagrange interpolation polynomial with nodes $x_1, \ldots,$ $x_{n+1}$. Assume that $|x_1| = u(\eta)$ and, for $\varepsilon > 0$, define $\tilde{\eta}$ as the measure with weights

$$(A.12) \qquad \tilde{w}_i = \frac{1}{2} \int_{-1}^1 \prod_{j=1, j \neq i}^{n+1} \frac{x - \tilde{x}_j}{\tilde{x}_i - \tilde{x}_j} \, dx, \qquad i = 1, \ldots, n+1,$$

at the points $\tilde{x}_1 = x_1 + \varepsilon, \tilde{x}_2 = x_2, \ldots, \tilde{x}_{n+1} = x_{n+1}$. If $\varepsilon$ is sufficiently small, it follows that all weights $\tilde{w}_j$ are positive, which implies $\tilde{\eta} \in \Xi_{n+1}$ and $\tilde{x}_1 > x_1$. In the case $|x_{n+1}| = u(\eta)$ we apply exactly the same argument to the point $x_{n+1}$ and obtain a measure $\tilde{\eta}$ with $u(\tilde{\eta}) < u(\eta)$. Consequently, the infimum in (A.11) cannot be attained in $\Xi_{n+1}$, that is,

$$(A.13) \qquad \inf_{\eta \in \Xi} u(\eta) = \inf_{\eta \in \Xi_n} u(\eta).$$

We now prove that the infimum on the right-hand side cannot be attained in $\Xi_t$ whenever $t > \lfloor n/2 \rfloor + 1 = r$. For this, consider a measure

$$\eta = \begin{pmatrix} x_1 & \cdots & x_t \\ w_1 & \cdots & w_t \end{pmatrix} \in \Xi_t.$$

Recall from Lemma 3.1 that the weights have to satisfy (3.4) (with $t = r$) and that by (3.3) the support points satisfy

$$(A.14) \qquad \int_{-1}^1 \prod_{j=1}^t (x - x_j) x^\ell \, dx = 0, \qquad \ell = 0, \ldots, n - t.$$

Without loss of generality we assume that $|x_1| = u(\eta)$ and we will again construct a design with a smaller value of $u(\eta)$. For this we define $\tilde{x}_1 = x_1 + \varepsilon$, $\tilde{x}_j = x_j$ ($j = 3 + n - t, \ldots, t$),

$$(A.15) \qquad \psi(x) = (x - \tilde{x}_1) \prod_{j=3+n-t}^t (x - \tilde{x}_j)$$

and construct the remaining points, say $\tilde{x}_2, \ldots, \tilde{x}_{n-t+2}$, such that the resulting measure with weights of the form (3.4) at the points $\tilde{x}_1, \ldots, \tilde{x}_t$ defines a quadrature formula of degree $n$. A necessary condition for this property is given by

$$(A.16) \qquad \int_{-1}^1 \psi(x) \prod_{j=2}^{n-t+2} (x - \tilde{x}_j) x^\ell \, dx = 0, \qquad \ell = 0, \ldots, n - t,$$

and a straightforward calculation shows that these equations are equivalent to the system

$$(A.17) \qquad \frac{\partial}{\partial \tilde{x}_i} \int_{-1}^1 \psi(x) \prod_{j=2}^{n-t+2} (x - \tilde{x}_j)^2 \, dx = 0, \qquad i = 2, \ldots, n-t+2,$$

which determines the points $\tilde{x}_j = \tilde{x}_j(\varepsilon)$ $(j = 2, \ldots, n-t+2)$ implicitly as a function of the parameter $\varepsilon$, where $\tilde{x}_j(0) = x_j$ $(j = 2, \ldots, n-t+2)$. We now consider the Jacobi matrix of the system (A.17),

$$(A.18) \qquad J(\tilde{x}, \varepsilon) = \left( \frac{\partial^2}{\partial \tilde{x}_j \, \partial \tilde{x}_j} \int_{-1}^1 \psi(x) \prod_{k=2}^{n-t+2} (x - \tilde{x}_k)^2 \, dx \right)_{i,j=2}^{n-t+2},$$

where we use the notation $\tilde{x} = (\tilde{x}_2, \ldots, \tilde{x}_{n-t+2})$. This matrix can be calculated using the quadrature formula corresponding to the measure $\eta$. For this, note that $\tilde{x}(0) = (x_2, \ldots, x_{n-t+2})$ and define the polynomial

$$(A.19) \qquad g_i(x) = 2\psi(x) \prod_{j=2, j \neq i}^{n-t+2} (x - x_j)^2 \, dx.$$

Note that the degree of $g_i$ is less than or equal to $n$. If $e_i \in \mathbb{R}^i$ denotes the $i$th unit vector, we obtain

$$e_i^T J(\tilde{x}(0), 0) e_i = \int_{-1}^1 g_i(x) \, dx = \sum_{j=1}^t g_i(x_j) w_j = g_i(x_i) w_i$$

$$= 2 w_i (x_i - x_1) \prod_{j=n-t+3}^t (x_i - x_j) \prod_{j=2, j \neq i}^{n-t+2} (x_j - x_i)^2 \neq 0,$$

and by a similar calculation, $e_i^T J(\tilde{x}(0), 0) e_j = 0$ whenever $i \neq j$ $(i, j = 2, \ldots, n-t+2)$. Consequently, $J(\tilde{x}(0), 0)$ is diagonal with $\det J(\tilde{x}(0), 0) \neq 0$. It now follows from the implicit function theorem (see [9]) that for sufficiently small $\varepsilon$ there exist analytic functions $\tilde{x}_2(\varepsilon), \ldots, \tilde{x}_{n-t+2}(\varepsilon)$ such that (A.16) is satisfied for the points $\tilde{x}_1, \ldots, \tilde{x}_t$. This implies that for sufficiently small $\varepsilon$, the weights defined by (A.12) are positive and the design

$$(A.20) \qquad \tilde{\eta} = \begin{pmatrix} \tilde{x}_1 & \cdots & \tilde{x}_t \\ \tilde{w}_1 & \cdots & \tilde{w}_t \end{pmatrix}$$

defines a quadrature formula of degree $n$ with positive weights, that is, $\tilde{\eta} \in \Xi_t$. If necessary [in the case $|x_t| = u(\eta)$] we use a similar argument for the largest support point of the probability measure $\eta$ and finally obtain a probability measure $\tilde{\eta} \in \Xi_t$ such that $u(\tilde{\eta}) < u(\eta)$. In other words, the infimum on the right-hand side of (A.13) cannot be obtained in $\Xi_t$ if $t > \lfloor n/2 \rfloor + 1$. From (A.15) it is easy to see that this construction is not possible if $t < r = \lfloor n/2 \rfloor + 1$. Moreover, it is well known (see, e.g., [8]) that a quadrature formula that integrates polynomials of degree $n$ exactly must have at least $r$ support points. Moreover, if it has $r$ support points, it is uniquely determined and given by the probability measure $\eta^*$ defined in (A.8), which completes the proof of Lemma A.1. $\square$

**Acknowledgments.** Parts of this paper were written during Dette's visit to the Institute of Statistics in Louvain-la-Neuve and he would like to thank the Institute for its hospitality. The authors are also grateful to Isolde Gottschlich, who typed parts of this paper with considerable technical expertise, and to Kay Pilz for some computational assistance. The authors would also like to thank the reviewers for their constructive comments on an earlier version of this manuscript.

H. DETTE
FAKULTÄT FÜR MATHEMATIK
RUHR-UNIVERSITÄT BOCHUM
44780 BOCHUM
GERMANY
E-MAIL: holger.dette@ruhr-uni-bochum.de

V. B. MELAS
A. PEPELYSHEV
DEPARTMENT OF MATHEMATICS
ST. PETERSBURG STATE UNIVERSITY
ST. PETERSBURG
RUSSIA
E-MAIL: v.melas@pobox.spbu.ru
        andrey@ap7236.spb.edu